\documentclass[12pt,psamsfonts]{amsart}
\usepackage{amsmath}
\usepackage{amsthm}
\usepackage{amssymb}
\usepackage{amscd}
\usepackage{amsfonts}
\usepackage{amsbsy}

\def\a{\alpha}

\def\t{\theta}

\def\G{\Gamma}

\def\ve{{\varepsilon}}
\def\e{{\varepsilon}}
\def\Cs{\mbox{Cs\,}}
\def\Sn{\mbox{Sn\,}}
\def\ov{\overline}

\newcommand{\bbox}{\ \hfill\rule[-1mm]{2mm}{3.2mm}}

\newcommand{\R}{\ensuremath{\mathbb{R}}}

\newtheorem {theorem} {Theorem}
\newtheorem {proposition} [theorem]{Proposition}

\newtheorem {remark} [theorem]{Remark}

\begin{document}

\title[Nilpotent and degenerate centers]
{The problem of distinguishing between a center and a focus for
nilpotent and degenerate analytic systems}

\author[H. Giacomini, J. Gin\'{e} and J. Llibre]
{Hector Giacomini$^1$, Jaume Gin\'e$^2$ and Jaume Llibre$^3$}

\address{$^1$  Laboratoire de Math\'ematique et Physique Th\'eorique,
CNRS (UMR 6083), Facult\'e des Sciences et Techniques,
Universit\'e de Tours, Parc de Grandmont, 37200 Tours, France}

\email{giacomini@phys.univ-tours.fr}

\address{$^2$ Departament de Matem\`atica, Universitat de Lleida,
Av. Jaume II, 69. 25001 Lleida, Spain}

\email{gine@eps.udl.es}

\address{$^3$ Departament de Matem\`{a}tiques,
Universitat Aut\`{o}noma de Barcelona, 08193 Bellaterra, Barcelona,
Spain}

\email{jllibre@mat.uab.es}

\subjclass{Primary 34C05. Secondary 58F14.}

\keywords{nilpotent center, degenerate center, Liapunov constants,
cyclicity}
\date{}
\dedicatory{We dedicate this paper to the memory of Javier
Chavarriga}

\maketitle

\begin{abstract}
In this work we study the centers of planar analytic vector fields
which are limit of linear type centers. It is proved that all the
nilpotent centers are limit of linear type centers and
consequently the Poincar\'e--Liapunov method to find linear type
centers can be also used to find the nilpotent centers. Moreover,
we show that the degenerate centers which are limit of linear type
centers are also detectable with the Poincar\'e--Liapunov method.
\end{abstract}

\section{Introduction and statement of the main results}\label{s1}

Two of the main and oldest problems in the qualitative theory of
differential systems in $\R^2$ is the distinction between a center
and a focus, called the {\it center problem}; and the
determination of the first integrals in the case of centers, see
for instance \cite{Anosov}. This paper deals with these two
problems for the class of analytic differential systems.

Let $p\in \R^2$ be a singular point of a differential system in
$\R^2$. We say that $p$ is a {\it center} if there is a
neighborhood $U$ of $p$ such that all the orbits of $U\setminus
\{p\}$ are periodic, and we say that $p$ is a {\it focus} if there
is a neighborhood $U$ of $p$ such that all the orbits of
$U\setminus \{p\}$ spiral either in forward or in backward time to
$p$.

Once we have a center at $p$ of a differential system in $\R^2$,
another problem is to know if there exists or not a first integral
$H$ defined in some neighborhood $U$ of $p$ (i.e. a non--constant
function $H:U\to \R$ such that $H$ is constant on the orbits of
the differential system), and to know the differentiability of $H$
with respect to the differentiability of the system. More
specifically, we assume that we have an analytic differential
system having a center at $p$. Then, it is known that there exists
a $C^{\infty}$ first integral defined in some neighborhood of $p$,
see \cite{Sa}. It is also known that there exists an analytic
first integral defined in $U\setminus \{p\}$ for some neighborhood
$U$ of $p$, see \cite{LLNZ}; but such analytic first integral in
general cannot be extended to $p$. For any center $p$ of an
analytic differential system in $\R^2$ it is an open problem to
characterize when there exists an analytic first integral in a
neighborhood of $p$, or simply a {\it local analytic first
integral at $p$}.

A singular point $p$ is a {\it monodromy} singular point of a real
analytic differential system in $\R^2$ if there is no {\it
characteristic orbit} associated to it; i.e., there is no orbit
tending to the singular point with definite tangent at this point.
Let $p$ be a singular point of an analytic differential system. If
$p$ is monodromy, then it is either a center or a focus, see
\cite{Ecalle, Ilyashenko}. Moreover, $p$ is a center if and only
if there exists a $C^{\infty}$ first integral defined in some
neighborhood of $p$, see \cite{Sa}.

Let $p\in \R^2$ be a singular point of an analytic differential
system in $\R^2$, and assume that $p$ is a center. Without loss of
generality we can assume that $p$ is the origin of coordinates (if
necessary we do a translation of coordinates sending $p$ at the
origin). Then, after a linear change of variables and a rescaling
of the time variable (if necessary), the system can be written in
one of the following three forms:
\begin{eqnarray}
\dot{x}&=&-y+ F_1(x,y), \quad \, \, \dot{y}=x+F_2(x,y); \label{l1}\\
\dot{x}&=&y+ F_1(x,y), \quad \quad \, \dot{y}=F_2(x,y); \label{ll1}\\
\dot{x}&=&F_1(x,y), \quad \quad \quad \quad \dot{y}=F_2(x,y);
\label{ll2}
\end{eqnarray}
where $F_1(x,y)$ and $F_2(x,y)$ are real analytic functions
without constant and linear terms, defined in a neighborhood of
the origin. In what follows a center of an analytic differential
system in $\R^2$ is called {\it linear type, nilpotent or
degenerate} if after an affine change of variables and a rescaling
of the time it can be written as system (\ref{l1}), (\ref{ll1}) or
(\ref{ll2}), respectively.

The characterization of the linear type centers in terms of the
existence of an analytic first integral is due to Poincar\'e
\cite{Po} and Liapunov \cite{Ly}, see also Moussu \cite{Mos}.

\medskip

\noindent {\bf Linear Type Center Theorem.} {\it The real analytic
differential system \eqref{l1} has a center at the origin if and
only if there exists a local analytic first integral of the form
$H=x^2+y^2+ F(x,y)$ defined in a neighborhood of the origin, where
$F$ starts with terms of order higher than $2$.}

\medskip

An analytic system on the plane will have a singular point of
center type if a countable number of conditions on the
coefficients of the system are satisfied, see \cite{Anosov}. Based
on the Linear Type Center Theorem there is a method, called the
{\it Poincar\'{e}--Liapunov method}, which consists in determining
when a system of the form (\ref{l1}) has a local analytic first
integral at the origin, and consequently a center at this point.
This algorithm looks for a formal power series of the form
\begin{equation}\label{3}
H(x,y)=\sum_{n=2}^\infty H_n(x,y),
\end{equation}
where $H_2(x,y)=(x^2+y^2)/2$, and for each $n$, $H_n(x,y)$ are
homogeneous polynomials of degree $n$, so that
\begin{equation}\label{333}
\dot{H}=\sum_{k=2}^{\infty} V_{2k} (x^2+y^2)^k,
\end{equation}
where the $V_{2k}$'s are called the {\it Liapunov constants}. It
is known that the Liapunov constants are polynomials in the
coefficients of system (\ref{l1}). We note that the
Poincar\'{e}--Liapunov method for analytic differential systems is an
algorithm which at each step uses only a finite jet of the system
for the calculation of a Liapunov constant. The singular point is
a center if and only if all the Liapunov constants vanish. For
more details see \cite{Anosov} and references therein.

Until now there is no algorithm comparable to the
Poincar\'e--Liapunov method for determining the center conditions
in the case of nilpotent and degenerate singular points, except if
the singular point has no characteristic direction because in this
last case we can use the algorithm of Bautin \cite{Ba} (see also
\cite{CG1, Anosov, Mo}). In any case the necessary computations
for applying Bautin's algorithm are in general more difficult to
implement that the ones coming from the Poincar\'e--Liapunov
method. In this paper we shall show that essentially the
Poincar\'{e}--Liapunov algorithm also works for determining the
analytic nilpotent centers and a subclass of the analytic
degenerate centers.

Our main result is the following one.

\begin{theorem}[{\bf Nilpotent Center Theorem}]\label{t1}
Suppose that the origin of the real analytic differential system
\eqref{ll1} is a center, then there exist analytic functions $G_1$
and $G_2$ without constants terms, such that the system
\begin{equation}\label{PL}
\dot{x}=y+ F_1(x,y) + \e x G_1(x,y), \quad  \, \dot{y}= -\e
x+F_2(x,y) + \e x G_2(x,y),
\end{equation}
has a linear type center at the origin for all $\e>0$.
\end{theorem}

Roughly speaking Theorem \ref{t1} can be stated saying simply that
{\it an analytic nilpotent center is always limit of analytic
linear type centers}. Theorem \ref{t1} is proved in Section
\ref{s2}.

By the Linear Type Center Theorem, system \eqref{PL} has a local
analytic first integral $H_{\ve}(x,y)$ at the origin for $\ve>0$.
If there exists
\[
\lim_{\ve\searrow 0} H_{\ve}(x,y),
\]
and it is a function $H(x,y)$ well defined in a neighborhood of
the origin, then $H(x,y)$ is a local first integral of system
\eqref{ll1} at the origin. Note that, in general, $H$ is not
analytic, see Remark \ref{rrr}.

We note that the Nilpotent Center Theorem reduces the study of the
nilpotent centers to the case of linear type centers. So, we can
apply the Poincar\'{e}--Liapunov method to system \eqref{PL}, looking
for analytic first integrals of the form $H= (\e x^2+y^2)/2+
F(x,y,\e)$, where $F$ starts with terms of order higher than $2$
in the variables $x$ and $y$. We determine the Liapunov constants
$V_{2k}$ from  \eqref{333}. Several examples showing the
application of the Poincar\'{e}--Liapunov method to detect nilpotent
centers are given in Section \ref{s4}.

Based in the results obtained for nilpotent centers we establish
the following definition.

Suppose that the origin of the real analytic differential system
\eqref{ll2} is a center. We say that it is {\it limit of linear
type centers} if there exist $G_1$ and $G_2$ analytic functions in
$x$, $y$ and $ \e$, without constants and linear terms in $x$ and
$y$, such that the system
\begin{equation}\label{deff3}
\dot{x}= \e y+ F_1(x,y) + \e G_1(x,y,\e), \quad \dot{y}= -\e x+
F_2(x,y)+ \e G_2(x,y,\e),
\end{equation}
has a linear type center at the origin for all $\e \ne 0$
sufficiently small. A more general definition of limit of linear
type centers would be to consider functions $G_1$ and $G_2$ that
are not analytic in $\e$.

\begin{theorem}\label{pp3}
Suppose that the origin of the real analytic differential system
\eqref{ll1} or \eqref{ll2} is monodromy, and that this system is
limit of linear type centers of the form \eqref{PL} or
\eqref{deff3}, respectively. Suppose also that there are no
singular point of \eqref{PL} or \eqref{deff3} tending to the
origin when $\e$ tends to zero. Then, system \eqref{ll1} or
\eqref{ll2} has a center at the origin.
\end{theorem}

Theorem \ref{pp3} is proved in Section \ref{s2}. The condition
that there are no singular point tending to the origin when $\e$
tends to zero is easily verifiable using the lower order terms of
the perturbed system \eqref{PL} or \eqref{deff3}.

Another difficulty of the problem of distinguishing between a
center and a focus becomes from the fact that this problem for
degenerate centers can be no {\it algebraically solvable}; i.e.,
it does not exist an infinite sequence of independent polynomial
expressions involving the coefficients of the system, such that
their simultaneous vanishing guarantees the existence of a center,
see \cite{Ar2, Anosov, I, IY}.

The problem of distinguishing between a center and a focus is
algebraically solvable in the class of analytic differential
systems of type (\ref{l1}) and (\ref{ll1}), see \cite{Anosov, IY,
Ly, Po}. {\it The Nilpotent Center Theorem provides a new proof of
the fact that nilpotent analytic centers are algebraically
solvable}. This is due to the fact that we have seen in Theorem
\ref{t1} that the nilpotent centers can be $\e$--approximated by
systems having linear type centers. Then, applying the
Poincar\'e--Liapunov method to these linear type centers, and
doing the limit when $\varepsilon \searrow 0$ we obtain algebraic
conditions characterizing the existence of a nilpotent center. So,
in fact {\it Theorem \ref{t1} provides an algorithm for solving
the center problem for nilpotent centers}. See the end of Section
\ref{s2} for more details, and also Section \ref{s4}.

For centers of the form (\ref{ll2}) some preliminary results exist
for distinguishing between a center and a focus, see for instance
\cite{Anosov, GLMM, GMM, G, M}.

We say that an analytic differential system in the plane is {\it
time--reversible} (with respect to an axis of symmetry through the
origin) if after a rotation
$$
\begin{array}{rccl}
\left( \begin{array}{c} \xi \\ \eta \end{array} \right) & = &
\left( \begin{array}{cc} \cos\alpha & -\sin\alpha \\ \sin\alpha &
\cos\alpha
\end{array} \right) &
\left( \begin{array}{c} x \\ y \end{array} \right) ,
\end{array}
$$
the system in the new variables $(\xi,\eta)$ becomes invariant by
a transformation of the form $(\xi,\eta,t)\mapsto (\xi,-\eta,-t)$.
The phase portrait of this new system is symmetric with respect to
the straight line $\xi=0$. We note that for all reversible
nilpotent centers which are symmetric with respect to a straight
line through the origin, this line can be only the line of the
axes $x$ or $y$. We remark that all the nilpotent centers that we
know are time--reversible or have an analytic first integral at
the origin.

In the case of degenerate centers is much more difficult to
distinguish between a center and a focus than in the case of
linear and nilpotent type centers. In the next theorem we present
some results for the degenerate centers.

\begin{theorem} \label{t3}
For a degenerate analytic center the following statements hold.
\begin{itemize}
\item[(a)] A Hamiltonian degenerate center is always limit of
linear type Hamiltonian centers.
\item[(b)] A time--reversible degenerate center is always limit
of linear type time--reversible centers.
\item[(c)] There are degenerate centers which are neither
Hamiltonian nor time--reversible that are limit of linear type
centers.
\item[(d)] Non algebraically solvable degenerate centers are not
limit of linear type centers.
\item[(e)] There are algebraically solvable degenerate centers
which are not limit of linear type centers.
\item[(f)] There exist degenerate centers with characteristic
directions which are limit of degenerate centers without
characteristic directions.
\end{itemize}
\end{theorem}

Theorem \ref{t3} is proved in Section \ref{s3}.

Let \eqref{ll2} be a family of analytic systems depending on
several parameters. Inside the degenerate centers of this family
we can determine those which are limit of linear type centers of
the form \eqref{deff3}. For this kind of systems we can apply the
Poincar\'{e}--Liapunov method to system \eqref{deff3} with $\e \ne
0$ and compute their Liapunov constants. Vanishing these Liapunov
constants we obtain the center conditions for the system
\eqref{deff3}. Taking the limit when $\e \to 0$ in these
conditions, we get the center conditions for the degenerate
centers \eqref{ll2}. System (\ref{deff3}) must be a linear type
center only {\it for $\e \ne 0$ sufficiently small}. In
consequence, in the applications of the Poincar\'{e}--Liapunov
algorithm, it is sufficient to calculate the Liapunov constants up
to first order in $\e$. In contrast, for the nilpotent centers,
which are always limit of linear type centers, we can calculate up
to any order in $\e$ because system (\ref{PL}) has a center at the
origin {\it for all $\e>0$}, and from this fact we can obtain
several conditions at each step of the algorithm.

In Section \ref{s4} we provide an example of the application of
this method to a family of polynomial differential systems
\eqref{ll2}.

Finally, in Section \ref{s5} we obtain some results on the
cyclicity of nilpotent and degenerate centers which are limit of
linear type centers. In particular, we prove the following result
(for a definition of cyclicity of a center see Section \ref{s5}).

\begin{proposition}\label{pp1} Consider a nilpotent center or a
degenerate center of a polynomial differential system \eqref{ll1}
or \eqref{ll2} of degree $m$. We suppose that this center is limit
when $\e \to 0$ of linear type centers of polynomial differential
systems of degree $n$ of the form \eqref{PL} or \eqref{deff3},
respectively. If the Liapunov constants of a general perturbation
of the same degree $n$ of the linear type centers \eqref{PL} or
\eqref{deff3} are well--defined when $\ve\to 0$ and the Poincar\'e
map for a perturbation of the initial nilpotent or degenerate
center is analytic, then the following statements hold:
\begin{itemize}
\item[(a)] The cyclicity of the nilpotent center \eqref{ll1} is at
most the cyclicity of the linear type center \eqref{PL} for all
$\e>0$. \item[(b)] The cyclicity of the degenerate center
\eqref{ll2} is at most the cyclicity of the linear type centers
\eqref{deff3} for $\e \ne 0$ sufficiently small.
\end{itemize}
\end{proposition}

\section{Proof of Theorems \ref{t1} and \ref{pp3}}\label{s2}

The characterization of the nilpotent centers in terms of the
existence of a symmetry is due to Berthier and Moussu \cite{BM}
who obtained the following result. We shall need it in the proof
of Theorem \ref{t1}.

\begin{theorem}\label{t2}
If the analytic system \eqref{ll1} has a center at the origin,
then there exists an analytic change of variables such that the
new system has also the form \eqref{ll1} and it is invariant by
the change of variables $(x,y,t)\to (-x,y,-t)$.
\end{theorem}

We recall from \cite{CGGL} that if the analytic system \eqref{ll1}
has a center at the origin and there exists an analytic change of
variables such that the new system has also the form \eqref{ll1}
and it is invariant by the change of variables $(x,y,t)\to
(x,-y,-t)$, then the system has an analytic first integral defined
in a neighborhood of the origin.

\medskip

\noindent {\it Proof of Theorem} \ref{t1}: Assume that the origin
of system \eqref{ll1} is a center. Theorem \ref{t2} and its proof
says that for any nilpotent center \eqref{ll1} corresponding to an
analytic vector field $X(x,y)$, there exists an analytic change of
variables $(x,y)\to (u,v)$ of the form
\begin{equation}\label{xy}
x=u+\ldots, \qquad y=v+\ldots,
\end{equation}
such that $X(x,y)$ written in the new variables is a vector field
of the form
\begin{equation} \label{aa1}
Y(u,v)= (v + {\ov F}_1(u,v), {\ov F}_2(u,v)),
\end{equation}
where ${\ov F}_1$ and ${\ov F}_2$ are analytic functions starting
with terms of second degree in $x$ and $y$, and the associated
differential system is invariant under the change of variables
$(u,v,t)\mapsto (-u,v,-t)$.

Now we consider the following perturbation of the vector field
(\ref{aa1}):
\begin{equation} \label{aa}
Y_{\ve}(u,v)= (v + {\ov F}_1(u,v), -\e u+ {\ov F}_2(u,v)),
\end{equation}
with $\e>0$. Since the eigenvalues at the singular point located
at the origin are $\pm \sqrt{\ve}\, i$, and the differential
system associated to the vector field (\ref{aa}) is invariant
under the change of variables $(u,v,t)\mapsto (-u,v,-t)$ (because
the unperturbed system is invariant), it follows that the origin
of the vector field (\ref{aa}) is a linear type center for all
$\ve>0$.

Using the inverse of the change of variables (\ref{xy}) we get
that the differential system associated to the vector field
(\ref{aa}) becomes
\begin{equation} \label{aa2}
\dot{x}= y + F_1(x,y)+\e x G_1(x,y), \quad \dot{y}= -\e
x+F_2(x,y)+\e x G_2(x,y),
\end{equation}
where $G_1$ and $G_2$ are analytic functions without constants
terms, depending on the change of variable (\ref{xy}). Let
$X_\e(x,y)$ be the vector field associated to system \eqref{aa2}.
Since $Y_\e(u,v)$ has a linear type center at the origin for all
$\e>0$, the same holds for $X_\e(x,y)$. This completes the proof
of the theorem. \bbox

\medskip

\noindent {\it Proof of Theorem} \ref{pp3}: Consider an analytic
system $(P,Q)$ of the form \eqref{ll1} or \eqref{ll2} with a
monodromy singular point $p$ at the origin. Suppose that this
system is limit of linear type centers $(P_\e,Q_\e)$ of the form
\eqref{PL} or \eqref{deff3}, respectively. Since the origin is
monodromy, if $S$ is a sufficiently small curve with an endpoint
at the origin, then the Poincar\'e map $\Pi: S \to S$ associated
to the system $(P,Q)$ is well--defined and the leading term is
always linear for a suitable choice of a semi--transversal
algebraic curve, which can have a singularity at the singular
point, see \cite{Anosov, Me}. The Poincar\'e map $\Pi_\e:S\to S$
associated to the system $(P_\e,Q_\e)$ is the identity for all
$\e>0$ if the center is nilpotent, and for $\e$ sufficiently small
and $\e\neq 0$ if the center is degenerate. Therefore, by the
theorem on analytic dependence on initial conditions and
parameters, it follows that $\Pi= \lim_{\e \searrow 0} \Pi_\e$ if
the center is nilpotent, or  $\Pi= \lim_{\e \to 0} \Pi_\e$ if the
center is degenerate. Hence, we conclude that $\Pi$ is the
identity. So, the monodromy singular point $p$ of $(P,Q)$ is a
center. The condition that there are not singular points tending
to the origin when $\e$ tends to zero guarantees that the domain
of $\Pi_\e$ does not reduce to the origin when $\e$ tends to zero.
\bbox

\medskip

Theorems \ref{t1} and \ref{pp3} can be used to detect nilpotent
centers of analytic differential systems applying the algorithm of
Poincar\'e--Liapunov. In the particular case of polynomial systems
the method works as follows. We consider the system
\begin{equation} \label{aa3}
\dot{x}=y+F_1(x,y), \quad \dot{y}=F_2(x,y),
\end{equation}
where $F_1$ and $F_2$ are polynomials without constants and linear
terms containing a set of arbitrary parameters and such that the
origin is a monodromy singular point. We recall that using
Andreev's Theorem we can know when a nilpotent singular point is
or not monodromy, see \cite{AN}. For detecting the centers of
(\ref{aa3}), according with Theorem \ref{t1}, we consider the
perturbed system
\begin{equation} \label{aa4}
\dot{x}=y+F_1(x,y) + \e x G_1(x,y), \quad \dot{y}= - \e x +
F_2(x,y)+ \e x G_2(x,y),
\end{equation}
where $xG_1$ and $xG_2$ are analytic functions starting with
quadratic terms in $x$ and $y$. We apply now the
Poincar\'e--Liapunov algorithm to determine necessary conditions
to have a center at the origin for system (\ref{aa4}). In general,
these conditions will be satisfied by choosing conveniently the
coefficients of the analytic functions $G_1$ and $G_2$. When this
is not possible we must employ the parameters of the polynomial
system (\ref{aa3}). In this way we will obtain necessary
conditions for the existence of a center at the origin of system
(\ref{aa3}). The set of sufficient conditions of center for the
non--perturbed system (\ref{aa3}) will be obtained in a finite
number of steps, because the Hilbert's basis theorem guarantees
that this process is finite. Every time that we find a necessary
condition for the non--perturbed system (\ref{aa3}) we must to
study if the non-perturbed system (\ref{aa3}) already have a
center at the origin. As the number of steps is finite and for
determining each Poincar\'e--Liapunov constant of the perturbed
system (\ref{aa4}) we need only a finite jet, the necessary
perturbation to detect the center cases will be polynomial, i.e.,
the functions $G_1$ and $G_2$ will be polynomials. We note that,
under the assumptions of Theorem \ref{pp3}, it is not possible to
satisfy the center conditions of (\ref{aa4}) only with the
parameters of the perturbation because in that case using Theorem
\ref{pp3} the nilpotent polynomial system would have a center for
arbitrary values of the parameters of the family, which is a
contradiction if the initial system (\ref{aa3}) has not a center
at the origin.

\section{Proof of Theorem \ref{t3}}\label{s3}

In this section we shall work with an analytic degenerate center
(\ref{ll2}) defined in a neighborhood of the origin.

\medskip

\noindent{\it Proof of Theorem \ref{t3}(a)}:  Suppose that system
(\ref{ll2}) is Hamiltonian with Hamiltonian $H= H(x,y)$. The
system
\begin{equation}\label{ll3}
\dot{x}=- \ve y + F_1(x,y), \quad \dot{y}= \ve x+ F_2(x,y),
\end{equation}
is also a Hamiltonian system with the Hamiltonian first integral
$\ve(x^2+y^2)/2+ H(x,y)$. Consequently, system (\ref{ll3}) has a
linear type center at the origin for $\ve \ne 0$, and the initial
degenerate center (\ref{ll2}) is obtained taking in system
(\ref{ll3}) the limit when $\varepsilon \to 0$ . \bbox

\medskip

\noindent{\it Proof of Theorem \ref{t3}(b)}: Without loss of
generality, taking into account the definition of a
time--reversible system, we can assume that system (\ref{ll2}) is
invariant by the change of variables $(x,y,t)\mapsto (x,-y,-t)$.
Consider the perturbation of it given by a system of the form
(\ref{ll3}). Then, it is easy to see that system (\ref{ll3}) is
also invariant under the change of variables $(x,y,t)\mapsto
(x,-y,-t)$. Therefore, since the eigenvalues of the linear part at
the origin of system (\ref{ll3}) are $\pm\sqrt{|\ve|}\, i$, it has
a linear type center at the origin for $\ve \ne 0$. Again, the
initial degenerate center (\ref{ll2}) is obtained taking in system
(\ref{ll3}) the limit when $\varepsilon \to 0$ . \bbox

\medskip

\noindent{\it Proof of Theorem \ref{t3}(c)}: Consider the
following quartic polynomial differential system
\begin{equation}
\begin{array}{lll}
\dot{x}&=& (-y+y^2)(x^2+y^2), \\
\dot{y}&=& (x+2x^2) (x^2+y^2).
\end{array}
\label{LL5}
\end{equation}
It is easy to see that this system has a degenerate center at the
origin, because removing the common factor $x^2+y^2$ (doing a
change of the independent variable) we get a quadratic Hamiltonian
system having a center at the origin.

It is easy to check that system (\ref{LL5}) is not Hamiltonian,
and that it has the first integral
\begin{equation}\label{LL7}
H(x,y)= (x^2+y^2)/2 +2x^3/3-y^3/3.
\end{equation}
Now we claim that system (\ref{LL5}) is not time--reversible.
Suppose that it is time--reversible. Then, there exists a rotation
which pass the variables $(x,y)$ to the new variables $(u,v)$,
given by
\begin{equation}\label{jiji}
u= \cos \a\, x-\sin \a\, y, \qquad v=\sin\a\, x+ \cos\a\, y,
\end{equation}
which transforms the axis of symmetry into the line $u=0$. In the
new variables the system becomes $\dot u= P(u,v)$ and $\dot v=
Q(u,v)$. Since this system must be invariant by $(u,v,t)\to
(u,-v,-t)$, we must have
\[
P(u,v)= -P(u,-v), \qquad  Q(u,v)= Q(u,-v).
\]
These two equations are satisfied if and only if
\[
\cos\a \sin\a(2\cos\a-\sin\a)= 0, \qquad 2\sin^3\a- \cos^3\a=0.
\]
Since this system has no solution, the claim is proved.

Now, consider the following perturbation of system (\ref{LL5})
\begin{equation}
\begin{array}{lll}\label{LL6}
\dot{x}&=& (-y+y^2)(x^2+y^2-\varepsilon), \\
\dot{y}&=& (x+2x^2)(x^2+y^2-\varepsilon).
\end{array}
\end{equation}
System (\ref{LL6}) has also the first integral (\ref{LL7}).
Consequently, system (\ref{LL6}) has a center at the origin for
all $\e\in \R$. This center is of linear type if $\e\ne 0$. Hence,
doing the limit $\varepsilon \to 0$ in system (\ref{LL6}), we
obtain the initial system (\ref{LL5}) with a degenerate center at
the origin. \bbox

\medskip

We have seen in Theorem \ref{t1} that the nilpotent centers can be
$\e$--approximated by systems having linear type centers. Then,
applying the Poincar\'e--Liapunov method to these linear type
centers, and doing the limit when $\varepsilon \searrow 0$ we
obtain algebraic conditions characterizing the existence of a
nilpotent center. We shall see this more explicitly in Section
\ref{s4}. Now, we shall show that this does not occur for
degenerate centers which are not algebraically solvable (proving
statement (d) of Theorem \ref{t3}), and for some classes of
degenerate centers which are algebraically solvable (proving
statement (e) of Theorem \ref{t3}).

\medskip

\noindent{\it Proof of Theorem \ref{t3}(d)}: It is known that the
problem of determining the center conditions at the origin of the
system
\begin{eqnarray*}
\dot x &=& x p_2-y p_1+4x(x^2 +\mu y^2)p_1, \\
\dot y &=& x p_1+y p_2+4y(x^2 +\mu y^2)p_1,
\end{eqnarray*}
where $p_1=x^2+a_4 x y+a_5 y^2$ and $p_2=a_1x^2+ a_2xy+a_3 y^2$,
is not algebraically solvable for some specific values of the
parameters, see \cite{I}. If such a system is limit of linear type
centers, it would be algebraically solvable. So, the proof of
statement (d) of Theorem \ref{t3} follows. \bbox
\medskip

\noindent{\it Proof of Theorem \ref{t3}(e)}: Consider the
following cubic homogeneous system
\begin{equation}
\begin{array}{lll}
\dot{x}&=& P(x,y)= 12 \lambda x^3 - 9 x^2 y - 20 \lambda x y^2 -
25 y^3
+ 9 \mu y^3, \\
\dot{y}&=& Q(x,y)= 9 x^3 + 12 \lambda x^2 y + 25 x y^2 - 20
\lambda y^3,
\end{array}
\label{ss5}
\end{equation}
with the monodromy condition that $xQ(x,y)$ $-yP(x,y)$ has no real
factors. This system has a degenerate center at the origin if and
only if $\mu=0$, or $\lambda=0$. This follows checking the
conditions (i) and (ii) of the Appendix which characterize the
homogeneous systems having a center at the origin. Condition (i)
is satisfied by the monodromy condition, and condition (ii) is
$\displaystyle \int_0^{2\pi} \frac{f(\t)}{g(\t)} d\t = 0$; i.e.
\begin{eqnarray*}
&&\frac{-4 \pi \lambda}{\sqrt{9\mu-25} \sqrt{81\mu+64}
\sqrt{-17-\sqrt{64+81\mu}}\sqrt{-17+\sqrt{64+81\mu}}}\cdot \\
&&\left( \sqrt{-17-\sqrt{64+81\mu}}
(160-27\mu-5\sqrt{64+81\mu})\right.+\\
&&\left. \sqrt{-17+\sqrt{64+81\mu}}(-160+27\mu-5\sqrt{64+81\mu})
\right)=0.
\end{eqnarray*}
It is easy to check that condition (ii) holds if and only if
$\mu=0$ or $\lambda=0$.  Moreover, these centers are algebraically
solvable because condition (ii) is algebraic. In the case $\mu=0$
condition (i) is directly satisfied because
$xQ(x,y)-yP(x,y)=(x^2+y^2)(9x^2+25y^2)$. In the case $\lambda=0$,
the homogeneous polynomial $xQ(x,y)-yP(x,y)$ is
$9x^4+34x^2y^2+25y^4-9\mu y^4$, and condition (i) is satisfied if
and only if $\mu<25/9$.

We consider the following perturbation of system (\ref{ss5})
\begin{equation}
\begin{array}{ll}\label{ss6}
&\dot{x}= -\e y+ 12 \lambda x^3 - 9 x^2 y - 20 \lambda x y^2 -
25 y^3 + 9 \mu y^3 + \e G_1, \\
&\dot{y}= \ \ \e x + 9 x^3 + 12 \lambda x^2 y + 25 x y^2 - 20
\lambda y^3 +\e G_2,
\end{array}
\end{equation}
where $G_i=G_i(x,y,\e)$, for $i=1,2$, are analytic functions in
$x$, $y$ and $\e$, without constants and linear terms in $x$ and
$y$. Applying the Poincar\'e--Liapunov method to system
(\ref{ss6}) (see Section \ref{s4} for more details), we obtain
that the first Liapunov constant is
\[
V_1=-8 \lambda + \e \bar{V}_1,
\]
where $\bar{V}_1$ is the first Liapunov constant of the analytic
system
\[
\begin{array}{lll}
\dot{x}&=& -\e y + \e G_1(x,y,\e), \\
\dot{y}&=& \ \ \e x + \e G_2(x,y,\e).
\end{array}
\]
As we must vanish $V_1$ up to first order in $\e$ we obtain only
the condition $\lambda =0$. So, the case $\mu=0$ and $\lambda \in
\mathbb{R}$ cannot be detected as limit of linear type centers.
\bbox

\medskip

\noindent{\it Proof of Theorem \ref{t3}(f)}: We consider the
system
\begin{equation}\label{ll8}
\dot{x}=-a \, y^3 , \quad \dot{y}=b\, x^5,
\end{equation}
with $ab>0$. It is a Hamiltonian system with Hamiltonian
\[
H(x,y)= \frac{a y^4}{4}+\frac{b x^6}{6}.
\]
It is easy to check that system \eqref{ll8} is a
$(2,3)$--quasi--homogeneous system of weight degree $8$ satisfying
conditions (i) and (ii) for having a degenerate center at the
origin, see the Appendix. Another way to see that the origin is a
degenerate center is noting that the level curves of $H$ are
ovals. It is easy to check that system (\ref{ll8}) has only one
characteristic direction, given by $y=0$.

Now, consider the following perturbation of system (\ref{ll8}):
\begin{equation} \label{ll9}
\dot{x}= -a y^3, \quad  \dot{y}= \ve x^3 + b x^5 ,
\end{equation}
with $a\, \e>0$. This system is also Hamiltonian, with
\[
H_\e(x,y)= \frac{a y^4}{4}+\frac{\ve x^4}{4}+\frac{b x^6}{6}.
\]
System (\ref{ll9}) has a degenerate center at the origin, because
the origin is surrounded by ovals. This system has no
characteristic direction. Now, doing the limit $\varepsilon \to 0$
in system (\ref{ll9}), we obtain the initial system (\ref{ll8})
with a degenerate center at the origin and with a
characteristic direction. \bbox \\

The example given in the proof of Theorem \ref{t3}(f) shows that,
in a similar way that we can apply the Poincar\'e--Liapunov method
to detect nilpotent centers, in the study of certain degenerate
centers with characteristic directions we can apply the Bautin
method for degenerate centers without characteristic directions
(see for instance \cite{Ba, Br}) to a convenient perturbation of
the system with characteristic directions.

\section{The Poincar\'e--Liapunov method for nilpotent
and degenerate systems}\label{s4}

In this section we illustrate how to apply the
Poincar\'{e}--Liapunov method to several families of polynomial
differential systems for detecting nilpotent or degenerate
centers. Some of these families have been studied recently by
other authors with different and more complicated techniques.
First we start studying some nilpotent centers. We note that the
simplest nilpotent polynomial centers must be of degree $3$,
because there are no nilpotent center for quadratic polynomial
differential systems, see for instance \cite{Sc2}. We consider the
system
\begin{equation}\label{g0}
\dot{x}= y + x^2 + k_2 x y,\qquad \dot{y}= k_1 x^2-x^3.
\end{equation}
We apply to this family the general algorithm, with a general
perturbation and we obtain the following result:

\begin{proposition}\label{CG01}
System \eqref{g0} has a nilpotent center at the origin if and only
if $k_1=k_2=0$.
\end{proposition}

\noindent{\it Proof}: Applying the Poincar\'e--Liapunov method to
the perturbed system
\begin{equation}\label{g00}
\dot{x}=  y + x^2 + k_2 x y + \e x G_1(x,y),\quad \dot{y}= -\e x
+k_1 x^2-x^3 + \e xG_2(x,y),
\end{equation}
where
\[
G_1(x,y) = \sum_{i+j \, \ge 1}^{\infty} a_{ij} x^i y^j, \quad
G_2(x,y) = \sum_{i+j \, \ge 1}^{\infty} b_{ij} x^i y^j,
\]
and $\e>0$, we obtain the first Liapunov constant
\[
V_1=\frac{2}{3+2\e+3\e^2} \left [ 2k_1+(2b_{10} +2a_{10}
k_1+b_{01} k_1-k_2)\e \right.
\]
\[ \left.
-(a_{01} -3a_{20} -2a_{10} b_{10} -b_{01} b_{10} -b_{11}+a_{10}k_2
)\e^2+(a_{02} -a_{01} a_{10} )\e^3 \right].
\]
We note that in $V_1$ only appear the linear and quadratic terms
of $G_1$ and $G_2$. Vanishing $V_1$ at any order in $\e$ we get
the necessary center condition $k_1=0$ for system \eqref{g0} and
for an arbitrary perturbation. We obtain also the conditions
$b_{10}=k_2/2$, $a_{01}=(6a_{20}+2b_{11}+b_{01}k_2)/2$, and
$a_{02}=(6a_{10}a_{20}-2a_{10}b_{11}+a_{10}b_{01}k_2)/2$ on the
parameters of the perturbation. The next Liapunov constant has the
form
\[
V_2=\frac{1}{(1+\e)(5-2\e+5\e^2)} \left [ -72 k_2 + O(\e) \right
].
\]
In the expression of $V_2$ we have contributions of the linear,
quadratic, cubic, quartic and quintic terms of $G_1$ and $G_2$.
Therefore, the conditions $k_1=k_2=0$ are necessary in order that
the origin of system \eqref{g0} be a center. These conditions are
also sufficient as it is explained in Remark \ref{rrr}. \bbox

We see that it has been sufficient to employ a polynomial
perturbation of degree $5$ ir order to determine the necessary and
sufficient conditions of center for system \eqref{g0}.

Although in Theorem \ref{t1} the perturbation is unknown, it is
surprising that with the simple perturbation $-\e x$ in $\dot{y}$
it is possible to obtain the center cases of many families, as it
will be shown in the following examples. We consider the system
\begin{equation}\label{g1}
\dot{x}= y + Axy + By^2,\qquad \dot{y}= -x^3+Kxy+Ly^3.
\end{equation}
Applying the Andreev results \cite{AN} we can see that the origin
of system \eqref{g1} is monodromy.

\begin{proposition}\label{CG1}
System \eqref{g1} has a nilpotent center at the origin if and only
if $AB-3L=0$ and $AB(A^2-2K)=0$.
\end{proposition}

\noindent{\it Proof}: Applying the Poincar\'e--Liapunov method to
the perturbed system
\begin{equation}\label{g11}
\dot{x}= y + Axy + By^2,\qquad \dot{y}= -\e x -x^3+Kxy+Ly^3,
\end{equation}
with $\e>0$, we obtain the first Liapunov constant
\[
V_1=-\frac{2\e^2(AB-3L)}{3+2\e+3\e^2}.
\]
Vanishing $V_1$ we get the first center condition $L=AB/3$. Now,
we compute the second Liapunov constant
\[
V_2= -\frac{2\e^2AB(A^2-2K)}{3(1+\e)(5-2\e+5\e^2)}.
\]
Vanishing $V_2$ we obtain the second center condition
$AB(A^2-2K)=0$. So, these two conditions are necessary in order
that the origin of the perturbed system \eqref{g11} be a center.
These two conditions are not necessary, in principle, for system
\eqref{g1}, because we must investigate for others polynomials
perturbations of the form
\[
\dot{x}= y + Axy + By^2 + \e x G_1(x,y),\ \ \dot{y}= -\e x
-x^3+Kxy+Ly^3 + \e x G_2(x,y).
\]
But, in \cite{CG1} it is proved that these two conditions are
necessary in order that the origin of system \eqref{g1} be a
center. We remark that in this particular system we do not need to
take $\e=0$ in the center conditions because they are independent
of $\e$.

Now we prove that these two conditions are sufficient. If $A=0$ or
$B=0$ (and consequently $L=0$), we have that system \eqref{g1} is
reversible with respect to $(x,y,t)\mapsto (-x,y,-t)$ or
$(x,y,t)\mapsto (x,-y,-t)$, respectively. Therefore, since the
origin is monodromy, it is a center.

If $AB \ne 0 $, $L=AB/3$ and $A^2-2K=0$, then the system has the
analytic first integral
\[
H=\exp(-Ax)\left( y^2- \frac{12}{A^4}- \frac{12}{A^3}x-
\frac{6}{A^2}x^2-  \frac{2}{A}x^3+ Axy^2+\frac23 By^3 \right).
\]
This first integral can be obtained using the theory of
integrability of Darboux, see for instance \cite{Ll}. In fact,
this first integral already appeared in \cite{CG1}. Since the
origin is monodromy, by the existence of this analytic first
integral defined at the origin it follows that the origin is a
center. We note that in this case the nilpotent center is neither
time--reversible nor Hamiltonian. \bbox

\medskip

We consider the system
\begin{equation}\label{g2}
\dot{x}= -y,\qquad \dot{y}= x^5+ax^6+y(bx^3+cx^4).
\end{equation}
Applying Andreev's results \cite{AN} we can see that the origin of
system \eqref{g2} is monodromy.

\begin{proposition}\label{CG2}
System \eqref{g2} has a nilpotent center at the origin if and only
if $ab=0$ and $c=0$.
\end{proposition}

\noindent{\it Proof}: Applying the Poincar\'e--Liapunov method to
the perturbed system
\begin{equation}\label{g22}
\dot{x}= -y,\qquad \dot{y}= \e x+ x^5+ax^6+y(bx^3+cx^4),
\end{equation}
with $\e>0$, we obtain the first Liapunov constant
\[
V_1= \frac{2\, \e\, c}{5+3\, \e+3\, \e^2+5\, \e^3}.
\]
Vanishing $V_1$ we get the first center condition $c=0$. Now, we
compute the second Liapunov constant
\[
V_2= -\frac{(2+7\, \e)\,a\,b}{128\, \e^2}.
\]
Vanishing $V_2$ we obtain the second center condition $ab=0$. So,
these two conditions are necessary in order that the origin of the
perturbed system \eqref{g22} be a center. These two conditions are
not necessary, in principle, for system \eqref{g2} because we must
investigate, as in the previous example, for others polynomials
perturbations of the form
\[
\dot{x}= -y + \e x G_1(x,y),\quad \dot{y}= \e x+
x^5+ax^6+y(bx^3+cx^4)+\e x G_2(x,y).
\]
But, in \cite{CG2} it is proved that these two conditions are
necessary in order that the origin of system \eqref{g2} be a
center.

Now we prove that these two conditions are sufficient. If $a=c=0$
or $b=c=0$, we have that the system is reversible with respect to
$(x,y,t)\mapsto (-x,y,-t)$ or $(x,y,t)\mapsto (x,-y,-t)$,
respectively. Therefore, since the origin is monodromy, it is a
center. \bbox

\medskip

Proposition \ref{CG1} is proved in \cite{CG1} by using Liapunov
polar coordinates, see \cite{Ly}, and computing some generalized
Liapunov constants. The method developed in \cite{CG1} is not
useful to solve the center problem of Proposition \ref{CG2}, see
\cite{CG1}. Proposition \ref{CG2} is proved in \cite{CG2} by using
the normal form theory and taking into account that a convenient
truncated normal form of the nilpotent system is a Lienard system.
The method developed in this paper solves both problems in a
unified form and in a more simple way, by computing the
Poincar\'e--Liapunov constants of a linear center type system. In
both proofs we have used the results of \cite{CG1} and \cite{CG2}
to prove that the conditions are necessary. This is not a
restriction of our method because we can apply it with a general
perturbation. But, in that case it is necessary to make a big
amount of computations for obtaining the necessary conditions.
This is the usual amount of computations that appear in the
application of the Poincar\'e--Liapunov method when the system
under study has several parameters.

\medskip

We consider the system
\begin{equation}\label{g3}
\dot{x}=
-y+a_{11}xy+a_{02}y^2+a_{30}x^3+a_{21}x^2y+a_{12}xy^2+a_{03}y^3
,\quad \dot{y}= x^3.
\end{equation}

\begin{proposition}\label{CG3}
System \eqref{g3} has a nilpotent center at the origin if and only
if $a_{30}=0$, $a_{02}a_{11}+a_{12}=0$, $a_{02}a_{11}a_{21}=0$,
and $a_{02}a_{11}a_{03}=0$.
\end{proposition}

We consider the system
\begin{equation}\label{g4}
\dot{x}= -y ,\quad \dot{y}=
a_{11}xy+a_{02}y^2+a_{30}x^3+a_{21}x^2y+a_{12}xy^2+a_{03}y^3.
\end{equation}

\begin{proposition}\label{CG4}
System \eqref{g4} has a nilpotent center at the origin if and only
if $a_{21}-a_{02}a_{11}=0$, $a_{03}=0$, $a_{02}a_{11} a_{30}=0$,
and $a_{02}a_{11}(3a_{02}^2 + 2a_{12})=0$.
\end{proposition}

The proofs of Propositions \ref{CG3} and \ref{CG4} are similar to
the proof of Propositions \ref{CG1} and  \ref{CG2} and we omit
them. Proposition \ref{CG3} is also proved in \cite{CG1} and
Proposition \ref{CG4} is proved in \cite{CG2}.

\medskip

As the previous examples show, in some cases, it is sufficient to
perturb system (\ref{ll1}) with $-\e x$ in $\dot{y}$, but there
are nilpotent centers which are limit of more general
perturbations and which cannot be detected only with the
perturbation $-\e x$ in $\dot{y}$, as we will see in the following
example.

We consider the system
\begin{equation}
\begin{array}{ll}
&\dot{x}= P(x,y)=y +xy +(1-a)y^2 + (1-a) xy^2-ax^4-ax^5, \\
&\dot{y}= Q(x,y)=cy^2-2x^3+cy^3-2x^3y+(c-2)x^4(1+y).
\end{array}
\label{cont}
\end{equation}
Applying the Andreev's Theorem \cite{AN} it is easy to see that
system (\ref{cont}) has a monodromy singular point at the origin.

\begin{proposition}\label{CG0}
System \eqref{cont} has a nilpotent center at the origin for all
values of $a$ and $c$.
\end{proposition}

\noindent{\it Proof}: System (\ref{cont}) has the following
analytic first integral
\[
H(x,y)=(1+x)^{-2c}(1+y)^{-2a}(x^4+y^2),
\]
which can be determined using the Darboux theory of integrability.
Therefore, system (\ref{cont}) has a center at the origin. \bbox

\medskip

Applying the Poincar\'e--Liapunov method to the perturbed system
\begin{equation}
\begin{array}{lll}
\dot{x}&=& P(x,y), \\
\dot{y}&=& -\e x+ Q(x,y).
\end{array}
\label{cont1}
\end{equation}
with $\e>0$, we obtain the first Liapunov constant
\[
V_1=\frac{2 \, \e^2c(1+2a)}{3+2\e+3\e^2}.
\]
Therefore, the first center condition is $c(1+2a)=0$. But, we know
that system (\ref{cont}) has a nilpotent center for all values of
$a$ and $c$. Then, it must exist another more general
$\e$--perturbation of system (\ref{cont}) which is a linear type
center for all values of $a$ and $c$. Consider the following
polynomial perturbed system
\begin{equation}
\begin{array}{lll}
\dot{x}&=& -\e x (ax+ax^2) +P(x,y) , \\
\dot{y}&=& -\e x (1+(1-c)x+y+(1-c)xy)+Q(x,y),
\end{array}
\label{cont2}
\end{equation}
where $P$ and $Q$ are defined in system (\ref{cont}). System
(\ref{cont2}) has a linear type center at the origin because it
has the following analytic first integral
\[
H(x,y)=(1+x)^{-2c}(1+y)^{-2a}(x^4+y^2+ \e x^2).
\]
Therefore, all the nilpotent centers of system (\ref{cont}) are
limit of the linear type centers of system (\ref{cont2}), but not
all the nilpotent centers of system (\ref{cont}) are limit of
linear type centers of system (\ref{cont1}). This example shows
that it is not always possible to obtain a nilpotent center as a
limit of linear type centers with the only perturbation $-\e x$ in
$\dot{y}$, even in the case where a local analytic first integral
exists.

\begin{remark}\label{rrr}
\rm{Consider the system
\begin{equation}\label{nil}
\dot{x} = y + x^2, \qquad \dot{y} = - \, x^3.
\end{equation}
Since this system is time--reversible with respect to the change
of variables $(x,y,t) \to (-x,y,-t)$, and the origin is monodromy
(see \cite{AN, Andronov}), it has a nilpotent center at the
origin. But it has neither a local analytic first integral, nor a
formal first integral defined at the origin, see the proof in
\cite{CGGL}.

Consider now the following perturbation of system (\ref{nil}).
\begin{equation}
\dot{x} = y + x^2, \qquad \dot{y} = -\ve x \, - x^3.
\label{ejm2.1}
\end{equation}
As this system is time--reversible with respect to the same change
of variables, for $\ve>0$ it has also a center at the origin.
Therefore, by the Linear Type Center Theorem we know that system
(\ref{ejm2.1}) has a local analytic first integral $H(x,y,\ve)$ at
the origin. It is possible to compute an explicit expression of it
given by
\[
\exp \left[2 \arg(\ve+x^2+i(x^2+2y-\ve)) \right](\ve^2 + x^4 - 2
\ve y + 2 x^2 y + 2 y^2).
\]
Now, taking the limit when $\ve \searrow 0$, we obtain the first
integral
\[
H(x,y)= \lim_{\ve \searrow 0} H(x,y,\ve)= \exp \left[2 \arg(
x^2+i(x^2+2y)) \right](x^4 + 2 x^2 y + 2 y^2),
\]
of system \eqref{nil}, which is not analytic at the origin.}  \qed
\end{remark}

We see in this example that the limit of an analytic first
integral defined in a neighborhood of the origin can be not
analytic.

In general the study of the nilpotent centers is easier with the
algorithm proposed in this work than applying the results of
\cite{BM}. In our case, we have two arbitrary functions $G_1$ and
$G_2$, while in the algorithm consequence of the results of
\cite{BM} there are three arbitrary functions, the one which
appear in the normal form for the nilpotent center and the two
coming from the change of variables. Moreover, for polynomial
systems, under the assumptions of Theorem \ref{pp3}, the two
arbitrary functions $G_1$ and $G_2$ of our method are always
polynomials and this fact does not happen in the algorithm based
on the results of \cite{BM}.

Now we apply the Poincar\'{e}--Liapunov method to detect degenerate
centers in a family of polynomial differential systems.

We consider the polynomial  system
\begin{equation}
\begin{array}{lll}
\dot{x}&=& -a(1 + x)(x^4 - 4y^3 - 3y^4) + \mu y^3, \\
\dot{y}&=&  -a(1 + y)(4x^3 + 3 x^4 - y^4) + \lambda x^5.
\end{array}
\label{g5}
\end{equation}
with the monodromy condition $a \mu >0$ if $a \ne 0$ and $\mu
\lambda <0$ if $a=0$.
\begin{proposition}\label{CG5}
System \eqref{g5}, with the monodromy condition $a \mu >0$ if $a
\ne 0$ and $\mu \lambda <0$ if $a=0$, has degenerate centers at
the origin which are limit of linear type centers of the form
\eqref{deff3} with $G_1=G_2=0$ if and only if $\mu=\lambda=0$, or
$a=0$.
\end{proposition}

\noindent{\it Proof}: Applying the Poincar\'e--Liapunov method to
the perturbed system
\begin{equation}
\begin{array}{lll}
\dot{x}&=& \e y -a(1 + x)(x^4 - 4y^3 - 3y^4) + \mu y^3, \\
\dot{y}&=& -\e x -a(1 + y)(4x^3 + 3 x^4 - y^4) + \lambda x^5,
\end{array}
\label{gg5}
\end{equation}
with $\e \ne 0$, we obtain the first Liapunov constant
\[
V_1= -\frac{ a \mu}{\e}.
\]
Vanishing $V_1$ we get the first center condition $a \mu=0$. Now,
we compute the second Liapunov constant
\[
V_2= -\frac{5 \, a \lambda}{8 \e}.
\]
Vanishing $V_2$ we obtain the second center condition $a
\lambda=0$. So, these two conditions are necessary in order that
the origin of the perturbed system (\ref{gg5}) be a center.
Therefore, these two conditions are necessary in order that the
origin of system \eqref{g5} be a center which is limit of linear
type centers of the form \eqref{deff3} with $G_1=G_2=0$.

Now we prove that these two conditions are sufficient. If $a=0$ we
have that the system is Hamiltonian and reversible with respect to
$(x,y,t)\mapsto (x,-y,-t)$. Therefore, since the origin is
monodromy (because it has no characteristic directions), it is a
center. Now, taking the limit when $\ve \to 0$, we obtain a
Hamiltonian system which has a degenerate center at the origin.

If $\mu=\lambda=0$ it can be shown that system \eqref{g5} has a
monodromy singular point at the origin. Moreover, it has the
analytic first integral
\[
H(x,y)=(1+x)^{-1}(1+y)^{-1}(x^4+y^4),
\]
defined in a neighborhood of the origin. Therefore, system
\eqref{g5} has a degenerate center at the origin. We note that
this degenerate center is neither time--reversible nor
Hamiltonian. \bbox

\medskip

By the same arguments that for the nilpotent centers of polynomial
differential systems, the degenerate centers of polynomial
differential systems which are limit of linear type centers of the
form \eqref{deff3} under the assumptions of Theorem \ref{pp3}, in
fact, are limit of linear type centers of the form \eqref{deff3}
where the two analytic functions $G_1$ and $G_2$ are always
polynomials.

\section{On the cyclicity of nilponent and degenerate centers}\label{s5}

Let $p$ be a center of a polynomial vector field of degree $m$.
The {\it cyclicity}, $c_n(p)$, of $p$ is the maximum number of
limit cycles, taking into account their multiplicity, that can
bifurcate from the singular point $p$ when we perturb it into the
class of all polynomial differential systems of degree $n\geq m$.

For a linear type center $p$ of a polynomial differential system
of degree $m$ it is known that if the number of its independent
Liapunov constants is $k$, then the cyclicity $c_m(p) \le k-1$ if
the Bautin ideal is radical, see for instance \cite{Ro}. Moreover,
if we perturb a polynomial differential linear type center of
degree $m$ and cyclicity $k-1$ inside the class of all polynomial
vector fields of degree $m$, we can get perturbed vector fields
with exactly $k-1$ hyperbolic limit cycles bifurcating from the
center. This is due to the relationship between the Liapunov
constants and the coefficients of the Poincar\'{e} map near a
center. For more details on this subject see \cite{Ro}.

As we have seen in the examples, in general, the Liapunov
constants are not well--defined when $\ve\to 0$, see for instance
the proof of Propositions \ref{CG2} and \ref{CG5}. Therefore, we
must impose that the limit of the Liapunov constants when $\ve\to
0$ be well--defined. If the Liapunov constants are well--defined
when $\ve\to 0$, then the Poincar\'e map obtained by the limit
$\Pi= \lim_{\e \to 0} \Pi_\e$ gives a formal series which can be
not convergent at any positive radius. Hence, we must also impose
to the Poincar\'e map to be convergent in a neighborhood of the
origin. Taking into account these conditions we can establish the
following result.

\begin{proposition}\label{p322}
Suppose that the origin of a polynomial differential system
\eqref{ll1} or \eqref{ll2} of degree $m$ is a center $p$, and that
this system is limit when $\e \to 0$ of polynomial differential
systems of degree $n$ of the form \eqref{PL} or \eqref{deff3},
respectively, which have linear type centers $p_{\e}$ at the
origin. If
\begin{itemize}
\item[(i)] the Liapunov constants of a general perturbation of the
same degree $n$ of the linear type centers \eqref{PL} or
\eqref{deff3} are well--defined when $\ve\to 0$, and \item[(ii)]
the limit of the Poincar\'e map when $\ve\to 0$ of the general
perturbation of the same degree $n$ of the linear type centers
\eqref{PL} or \eqref{deff3} is analytic in a neighborhood of the
origin,
\end{itemize}
then,
\begin{itemize}
\item[(a)] the cyclicity $c_n(p)$ of the nilpotent center \eqref{ll1} is at most

the cyclicity $c_n(p_{\e})$ of the linear type center \eqref{PL}
for all $\e>0$. \item[(b)] the cyclicity $c_n(p)$ of the
degenerate center \eqref{ll2} is at most the cyclicity
$c_n(p_{\e})$ of the linear type centers \eqref{deff3} for $\e \ne
0$ sufficiently small.
\end{itemize}
\end{proposition}

\noindent{\it Proof}: Let $(P,Q)$ be the polynomial vector field
of degree $m$ associated to the system of the form \eqref{ll1} or
\eqref{ll2} with a singular point $p$ of center type at the
origin. By assumptions, the vector field $(P,Q)$ is limit when $\e
\to 0$ of polynomial vector fields $(P_{\e},Q_{\e})$ of degree $n$
associated to systems of the form \eqref{PL} or \eqref{deff3},
respectively, which have linear type centers $p_{\e}$ at the
origin.

Let $(P_{\e}^*,Q_{\e}^*)$ be a general perturbed polynomial vector
field of degree $n$ of the vector field $(P_{\e},Q_{\e})$. Taking
the limit of $(P_{\e}^*,Q_{\e}^*)$ when $\e \to 0$, we obtain a
perturbed polynomial vector field $(P^*,Q^*)$ of degree at most
$n$ of the vector field $(P,Q)$. Since the Liapunov constants of
the system $(P_{\e}^*,Q_{\e}^*)$ are well--defined when $\ve\to
0$, and the Poincar\'e map $\Pi= \lim_{\e \to 0} \Pi_\e$  of
$(P^*,Q^*)$ is analytic in a neighborhood of the origin, we can
control the ciclicity of the polynomial vector field $(P^*,Q^*)$
by the Poincar\'e map $\Pi$ (with the same restrictions that for
linear type centers). Moreover, the number of independent Liapunov
constants of the system $(P^*,Q^*)$ is at most the number of the
independent Liapunov constants of the system
$(P_{\e}^*,Q_{\e}^*)$. Therefore, the cyclicity $c_n(p)$ of the
nilpotent center \eqref{ll1} is at most the cyclicity
$c_n(p_{\e})$ of the linear type center \eqref{PL} for all $\e>0$,
and the cyclicity $c_n(p)$ of the degenerate center \eqref{ll2} is
at most the cyclicity $c_n(p_{\e})$ of the linear type centers
\eqref{deff3} for $\e \ne 0$ sufficiently small. Then the
proposition follows. \bbox

\medskip

In general the degenerate problems present the more rich
structure. For instance, when we look for algebraic limit cycles
into the quadratic polynomial vector fields, there is one which is
given by a non--degenerate algebraic curve (the algebraic limit
cycle of degree 2), but there are many others (the algebraic limit
cycles of degree 4, 5, 6, and perhaps others) that are given by
degenerate algebraic curves. Here, a degenerate algebraic curve is
an algebraic curve having singular points. For more details about
algebraic limit cycles see \cite{CLSc}. However, in the
assumptions of Proposition \ref{p322}, it is clear that the
cyclicity of a non--linear type center $p$ which is limit of
linear type centers, is not more rich than the cyclicity of the
linear type centers. Hence, Proposition \ref{pp1} follows from
Proposition \ref{p322}.

\section{Appendix: Homogeneous and quasi--homogeneous systems}

In this appendix we introduce two classes of polynomial vector
fields having degenerate centers. For more details about them see
\cite{CL} and \cite{LLYZ}.

We consider {\it polynomial differential systems} in $\R^2$ of the
form
\begin{equation}\label{1}
\dot x = P(x,y), \qquad \dot y = Q(x,y),
\end{equation}
where $P$ and $Q$ are real polynomials in the variables $x$ and
$y$. We say that this system has {\it degree} $m$ if $m$ is the
maximum of the degrees of $P$ and $Q$.

If $P$ and $Q$ are coprime homogeneous polynomials of degree $m$,
then the centers of systems (\ref{1}) are characterized by: (i)
the homogeneous polynomial $xQ(x,y)-yP(x,y)$ has no real factors
(so $m$ is odd), and (ii)
\[
\int_0^{2\pi} \frac{f(\t)}{g(\t)} d\t= 0.
\]
Here
\begin{eqnarray*}
f(\t)&=& \cos\t P(\cos\t,\sin\t)+\sin\t Q(\cos\t,\sin\t), \\
g(\t)&=& \cos\t Q(\cos\t,\sin\t)- \sin\t P(\cos\t,\sin\t).
\end{eqnarray*}
Moreover, all the homogeneous centers are global centers; i.e. the
periodic orbits surrounding the center fulfill all $\R^2$ .

In what follows $p$ and $q$ always will denote positive integers.

We say that the function $H(x,y)$ is {\it
$(p,q)$--quasi--homogeneous of weight degree $m\geq 0$} if $H(\l^p
x,\l^q y)= \l^m H(x,y)$ for all $\l\in\R$.

We say that system (\ref{1}) is {\it $(p,q)$--quasi--homogeneous
of weight degree $m\geq 0$} if $P$ and $Q$ are
$(p,q)$--quasi--homogeneous functions of weight degrees $p-1+m$
and $q-1+m$, respectively. Note that the
$(1,1)$--quasi--homogeneous systems of weight degree $m$ are the
classical homogeneous polynomial differential systems of degree
$m$. We note that if system (\ref{1}) is
$(p,q)$--quasi--homogeneous, then the differential equation
$dy/dx= Q/P$ (another way to write system (\ref{1})) is invariant
by the change of variables $(x,y)\to (\l^p x,\l^q y)$.

If $P$ and $Q$ are coprime, then the centers of the
$(p,q)$--quasi--homogeneous systems (\ref{1}) of degree $m$ are
characterized by: (i) the $(p,q)$--quasi--homogeneous polynomial
$pxQ(x,y)- qyP(x,y)$ has no real factors, and (ii)
\[
\int_0^{2\pi} \frac{F(\t)}{G(\t)} d\t= 0.
\]
Here
\begin{eqnarray*}
F(\t)&=& \Cs^{2q-1}\t\, P(\Cs\t,\Sn\t)+\Sn^{2p-1}\t\, Q(\Cs\t,\Sn\t), \\
G(\t)&=& p\, \Cs\t\, Q(\Cs\t,\Sn\t)- q\,\Sn\t\, P(\Cs\t,\Sn\t),
\end{eqnarray*}
and $\Cs\t$ and $\Sn\t$ are the $(q,p)$--trigonometric functions.
Moreover, all the $(p,q)$--quasi--homogeneous centers are global
centers.

We recall that the {\it $(p,q)$--trigonometric functions} $z(\t)=
\Cs\t$ and $w(\t)= \Sn\t$ are the solution of the following
initial value problem
\[
\dot z= -w^{2p-1}, \quad \dot w= z^{2q-1}, \quad z(0)=
p^{-\frac{1}{2q}}, \quad w(0)= 0.
\]

It easy to check that the functions $\Cs\t$ and $\Sn\t$ satisfy
the equality
\[
p\, \Cs^{2q}\t+ q\, \Sn^{2p}\t= 1.
\]
For $p=q=1$ we have that $\Cs\t= \cos\t$ and $\Sn\t= \sin\t$; i.e.
the $(1,1)$--trigonometric functions are the classical ones. The
functions $\Cs\t$ and $\Sn\t$ are $\tau$--periodic functions with
\[
\tau= 2 \, p^{-\frac{1}{2q}}\,  q^{-\frac{1}{2p}}
\frac{\G(\frac{1}{2p}) \G(\frac{1}{2q})}
{\G(\frac{1}{2p}+\frac{1}{2q})}.
\]

\bigskip

\noindent{\bf Acknowledgments}.

The authors thank some suggestions of the referee that have
improved the present paper.

The second author is partially supported by a DGICYT grant number
MTM2005-06098-C02-02 and by a CICYT grant number 2005SGR 00550,
and by DURSI of Government of Catalonia ``Distinci\'o de la
Generalitat de Catalunya per a la promoci\'o de la recerca
universit\`aria". The third author is partially supported by a
DGICYT grant number MTM2005-06098-C02-01 and by a CICYT grant
number 2005SGR 00550.

\end{document}